\documentclass[12pt,english,twoside]{article}
\RequirePackage{amsmath,amssymb,amsfonts,textcomp,yhmath}
\RequirePackage{times} 
\RequirePackage{cases} 
\RequirePackage[pdftex]{hyperref}\hypersetup{backref,implicit=true,bookmarks=true,
colorlinks=true, linkcolor=blue,filecolor=blue,
pagecolor=blue,urlcolor=blue,pdfpagemode=UseOutlines}

\newtheorem{theorem}{Theorem}[section]
\newtheorem{proposition}[theorem]{Proposition}

\newtheorem{remark}[theorem]{Remark}
\newtheorem{prev}{Proof}[section]

\begin{document}
\markboth{\small{A. Nangue}}{\tiny{Einstein-Scalar Field System with
a cosmological constant on the type I Bianchi space-time}}
 \centerline{}


\begin{center}
\Large{{\bf Einstein-Scalar Field System with a cosmological
constant on the type I Bianchi  space-time}}
\end{center}

\centerline{}
 \centerline{\bf {Alexis Nangue}}
 \centerline{Department of Mathematics}
\centerline{Higher Teacher's Training College,}
  \centerline{University of Maroua, PO.Box 55, Maroua, Cameroon}
 \centerline{alexnanga02@yahoo.fr}

\begin{abstract}
In many cases a  scalar field can lead to accelerated expansion in
cosmological models. This paper contains mathematical results on
this subject particularly on type I Bianchi space-time. In this
paper, global existence to the coupled Einstein-scalar field system
which rules the dynamics of a kind of pure matter in the presence of
a scalar field and cosmological constant is proven.
\end{abstract}
\noindent
{\bf MR Subject Classification}: 83C05, 83F05 \\


\section{Introduction}
General Relativity is a theory of gravitation, which states that the
gravitational attraction that is observed between the masses caused
by deformation of space and time through these masses and not as an
attractive force between the masses as in the theory of Newton
gravitation law. From then space and time are no longer
indissociable. General Relativity abandons the notion of force and
replaces it with the concept of curvature of space-time. To be
complete, this theory must also provide a means of calculating the
curvature of space-time created by mass distribution. It does this
through a complex system of mathematical formulas : the Einstein
equations linking the geometry of space-time and properties of the
matter. In global dynamics, the search for solutions to Einstein
equations coupled in different material fields remains an active
area of research particularly. When considering the Einstein
equations on Space-time which have surface symmetry, we can always
eliminate any phenomenon of wave propagation by a  suitable
coordinates choice. But experiments have shown the existence of
gravitational waves. In the context of General Relativity,
gravitational waves are defined as disturbances of the metric that,
from the point view of the Einstein equations are decoupled from
disturbances of energy-momentum tensor. One way to model the
phenomenon of gravitational waves is to introduce a scalar field in
the gravitation sources. This is the focus of this paper. Several
studies have already been carried out on the notion of scalar field
like the works of \cite{ADR}, \cite{DT} \cite{HL} \cite{NNAN1} and
\cite {2}.
\\\indent
We choose the Einstein equations with constant cosmological; this
interest is a physical reason. Indeed, astrophysical observations,
based on the redshift light spectrum, showed that the universe was
accelerating expansion. It is the presence of the cosmological
constant in Einstein equations that can mathematically model this
phenomenon. An important part of General Relativity is cosmology,
which is the study of the structure and evolution of the whole
universe. The geometric frame selected here is the type I Bianchi
space-time of Generalizing the Robertson-Walker space-time is
homogeneous and isotropic : the latter being the  background area of
cosmology. The phenomena studied here are called homogeneous, that
is, they depend only on time. Indeed, in the space-time, observers
located on the same  constant time  hypersurface see exactly the
same events so that only the evolution over time to be really
significant. We study here the existence of a global solution, that
is, defined on,  $ [0; +\infty[$, of the homogeneous Einstein-scalar
field system on a type I Bianchi  space-time  with a perfect fluid
model pure radiation type.
\\\indent
Unless otherwise specified, Greek indices range from 0 to 3, and
Latin indices from 1 to 3. We adopt the Einstein summation
convention
$a_{\alpha}b^{\alpha}=\sum\limits_{\alpha=0}^{3}a_{\alpha}b^{\alpha}
$. We consider the Bianchi type I space-time $(\mathbb{R}^{4}, g)$
and we denote by $x^{\alpha}=(x^{0}, x^{i})=(t,x^{i})$, the usual
coordinates in $\mathbb{R}^{4}$; $g$  stands for the unknown metric
tensor of Lorentzian type with signature $(-, +, +, +)$ which can be
written :
\begin{equation}\label{1.1}
    g = - dt^{2} + a^{2}(t)(dx^{1})^{2} + b^{2}(t)\left[(dx^{2})^{2}+(dx^{3})^{2}\right]
\end{equation}
where $a>0 $ and $b>0$ are unknown functions of the single variable
$t$.\\\indent The Einstein-Scalar Field system with cosmological
constant reads as follows, according to \cite{1} :

\begin{numcases}\strut
R_{\alpha\beta}\; - \; \frac{1}{2}Rg_{\alpha\beta}\; +\; \Lambda
g_{\alpha\beta}\; = 8\pi (T_{\alpha\beta} + \tau_{\alpha\beta})
\label{1.2}\\
\nabla_{\alpha} \nabla^{\alpha}\phi=0 \label{1.3}\\
 T_{\alpha\beta} \; = \;\nabla_{\alpha}\phi\nabla_{\beta}\phi\; -
\;\frac{1}{2}g_{\alpha\beta}
\nabla^{\lambda}\phi\nabla_{\lambda}\phi
 \label{1.4}\\
 \tau_{\alpha\beta}\; = \;\frac{4}{3}\rho u_{\alpha}u_{\beta}\; +\;
\frac{1}{3}\rho g_{\alpha\beta}\label{1.5}
 \end{numcases}

where :
\begin{enumerate}
    \item[$\bullet$] (\ref{1.2}) are the Einstein equations for the
    metric tensor $g=(g_{\alpha\beta})$ which represents the
    gravitational field; $R_{\alpha\beta}$ is the Ricci tensor,
    contracted of the curvature tensor;
    $R=g^{\alpha\beta}R_{\alpha\beta}$ is the scalar curvature,
    contracted of the Ricci tensor.
    \item[$\bullet$](\ref{1.3}) is the wave equation in $\phi$ which
    represents the scalar field. Recall that $\nabla_{\alpha}$ is
    the covariant differentiation in $g$, and are raised and lowered
    following the rules : $V^{\alpha}=g^{\alpha\beta}V_{\beta}$; $V_{\alpha}=g_{\alpha\beta}V^{\beta}$,
    where $(g^{\alpha\beta}=g_{\alpha\beta})^{-1}$.
    \item[$\bullet$]The ordinary matter is modeling by (\ref{1.5}),
    which represents the relativistic perfect fluid of pure
    radiation type, in which $\rho\geq0$ is an unknown function of
    single variable $t$, representing the matter density. For simplicity,
     we consider a co-moving fluid, which means that
    $u^{i}=u_{i}=0$, where $u=(u^{\alpha})$ is a future time-like
    unit vector (i.e $g_{\alpha\beta}u^{\alpha}u^{\beta}=-1$, $u^{0}>0$).
    \item[$\bullet$](\ref{1.4}) represents the stress-matter-energy
    tensor associated to a scalar field $\phi$, which is as $\rho$ a
    real-valued function of $t$.
\end{enumerate}
\indent Now, recall that, solving the Einstein equations is
determining both the gravitational field and its sources : this
means that we have to determine every unknown function introduced
above, namely : $a$, $b$, $\rho$ and $\phi$. Notice that the
spatially homogeneous coupled Einstein-Scalar Field system turns out
to be a non linear second differential system. What we call global
solution in this paper, is a solution defined all over the interval
$[0, +\infty[$.

\vspace{5cm}

The paper is organized as follows :
\begin{enumerate}
    \item[$\bullet$] In section 2, we write the Einstein-Scalar
    field system in a explicit form.
    \item[$\bullet$] In section 3, we introduce the Cauchy problem
    and we prove the local existence of solutions.
    \item[$\bullet$] In section 4, we prove the global existence of
    existence.
\end{enumerate}

\section{Einstein-Scalar Field System in $a$, $b$, $\rho$, $\phi$}
In this section we are going to write the equations (\ref{1.2}) in
explicit form, and afterwards we proceed to a suitable change of
unknown functions. The evolution of solutions of the Einstein-Scalar
Field system with a cosmological constant on the Bianchi type I
space-time models, described by a perfect fluid with matter density
$\rho$, are governed following \cite{2}, by the constraint equation
\begin{equation}\label{1.6}
 \left(\dfrac{\dot{b}}{b}\right)^2+2
\dfrac{\dot{a}}{a}\dfrac{\dot{b}}{b}-\Lambda=8\pi\rho+4\pi\dot{\phi}^2,
\end{equation}
named Hamiltonian equation, \footnote{\tiny{Overdot denotes
differentiation with respect to time $t$.}} the evolution equations,
\begin{equation}\label{1.7}
\left(\dfrac{\dot{b}}{b}\right)^2+2
\dfrac{\ddot{b}}{b}-\Lambda=-\dfrac{8\pi\rho}{3}-4\pi\dot{\phi}^2,
\end{equation}
\begin{equation}\label{1.8}
\dfrac{\ddot{a}}{a}+
\dfrac{\dot{a}\dot{b}}{ab}+\dfrac{\ddot{b}}{b}-\Lambda=-\dfrac{8\pi\rho}{3}-4\pi\dot{\phi}^2,
\end{equation}
and the equations in $\phi$ and $\rho$, resulting from (\ref{1.3})
and conservation equation, given by :
\begin{equation}\label{1.9}
\ddot{\phi}\dot{\phi}+\left(\dfrac{\dot{a}}{a}+
2\dfrac{\dot{b}}{b}\right)\dot{\phi}^{2}=0
\end{equation}
and
\begin{equation}\label{1.10}
 \dot{\rho}+\dfrac{4}{3}\left(\dfrac{\dot{a}}{a}+ \dfrac{\dot{b}}{b}\right)\rho=
 0.
\end{equation}
In the next paragraphs, we study the local and global existence of
solutions $a$, $b$, $\rho$ and $\phi$ to the coupled system
(\ref{1.7}), (\ref{1.8}), (\ref{1.9}), (\ref{1.10}) subject to
constraint (\ref{1.6}). For this purpose, we make a change of
unknown functions in order to deduce an equivalent first order
differential system to which standard theory is applied. We set :
\begin{equation}\label{1.11}
    u=\frac{\dot{a}}{a}\;  ; \;v=\frac{\dot{b}}{b}\;  ;
    \;\psi=\frac{1}{2}\dot{\phi}^{2}.
\end{equation}
We deduce from (\ref{1.11}) :
\begin{equation}\label{1.12}
 \frac{\ddot{a}}{a}=u+u^2        \;  ; \; \frac{\ddot{b}}{b}=v+v^2.
\end{equation}
We choose to look for  a $\mathcal{C}^{2}$-non-decreasing scalar
field (i.e. $\dot{\phi}\geq 0$), then (\ref{1.11}) gives :
\begin{equation}\label{1.13}
    \dot{\phi}=\sqrt{2}\psi^{\frac{1}{2}}.
\end{equation}
According to (\ref{1.2}) and (\ref{1.3}) we deduce from (\ref{1.7}),
(\ref{1.8}), (\ref{1.9}), (\ref{1.10}) the equivalent first order
differential system :
\begin{numcases}\strut
\dfrac{du}{dt}=\dfrac{2}{3}\Lambda-u^{2}+\dfrac{1}{3}v^{2}-\dfrac{4}{3}uv-\dfrac{8}{3}\pi\psi
\label{1.14}\\
 \dfrac{d v}{dt}=\dfrac{2}{3}\Lambda-\dfrac{5}{3}v^{2}-\dfrac{1}{3}uv-\dfrac{8}{3}\pi\psi
 \label{1.15}\\
 \dfrac{d \phi}{dt}=\sqrt{2}\psi^{\frac{1}{2}}\label{1.16}\\
 \dfrac{d \rho}{dt}=-\dfrac{4}{3}(u+2v)\rho\label{17}\\
 \dfrac{d \psi}{dt}=-2(u+2v)\psi\label{18}
\end{numcases}
subject to the constraint :
\begin{equation}\label{19}
v^{2}+2uv-\Lambda=8\pi\rho+8\pi\psi,
\end{equation}
which we are going to study.
\section{Cauchy problem and constraint}
Let $a_{0}>0$, $\dot{a}_{0}$, $b_{0}>0$, $\dot{b}_{0}$, $\phi_{0}$
$\dot{\phi}_{0}>0$, $\rho_{0}$ be given real numbers. We look for
solutions $a$, $b$, $\rho$ and $\phi$ of the Einstein-Scalar field
system over $[0, T[$, $T\leq +\infty$ satisfying :
\begin{equation}\label{20}
    a(0)=a_{0}\; ; \;\dot{a}(0)=\dot{a}_{0}\; ; \;b(0)=b_{0}\; ; \;\dot{b}(0)=\dot{b}_{0}\; ; \;
    \phi(0)=\phi_{0}\; ; \;\dot{\phi}(0)=\dot{\phi}_{0}\; ; \;
    \rho(0)=\rho_{0}.
\end{equation}
Our objective now is to prove the local existence of solution
satisfying (\ref{20}), called  initial conditions, with The given
numbers $a_{0}$, $\dot{a}_{0}$, $b_{0}$, $\dot{b}_{0}$, $\phi_{0}$
$\dot{\phi}_{0}$, $\rho_{0}$   being the initial data.\\\indent It
is well known that equation (\ref{1.6}) called Hamiltonian
constraint is satisfied all over the domain of the solutions of
evolution equations, if and only if equation (\ref{1.6}) is
satisfied at $t=0$ i.e given (\ref{20}) if the initial data satisfy
:
\begin{equation*}
 \left(\dfrac{\dot{b}_{0}}{b_{0}}\right)^2+2
\dfrac{\dot{a}_{0}}{a_{0}}\dfrac{\dot{b}_{0}}{b_{0}}-\Lambda=8\pi\rho_{0}+4\pi\dot{\phi}_{0}^2,
\end{equation*}
which is calling the initial constraint. Now we are going to study
the equivalent first order differential system (\ref{1.14}) to
(\ref{18}), subject to constraint (\ref{18}) and with the initial
conditions at $t=0$, provided by (\ref{20}) :
\begin{equation}\label{21}
u(0):=u_{0}=\dfrac{\dot{a}_{0}}{a_{0}}; \;
v(0):=v_{0}=\dfrac{\dot{b}_{0}}{b_{0}};\,  \rho(0)=\rho_{0} \;
 ;\psi(0):=\psi_{0}=\dot{\phi}_{0};  \;  \phi(0)=\phi_{0}.
\end{equation}

\section{Local existence of solutions}
We use an iterative scheme.
\subsection{Construction of the iterated sequence}
We construct the sequence $S_{n}=(u_{n}, v_{n}, \rho_{n}, \psi_{n},
\phi_{n})$, $n\in \mathbb{N} $, as follows :
\begin{enumerate}
    \item[$\bullet$] Set $u_{0}=u(0); \;
v_{0}=v(0);\,  \rho_{0}=\rho(0) \;
 ;=\psi_{0}=\psi(0);  \;  \phi_{0}=\phi(0)$ where  $u_{0}$,
$v_{0}$,  $\rho_{0}$,
 $\psi_{0}$,   $\phi_{0}$
 are initial data which satisfy constraint equation.
    \item[$\bullet$]Define $S_{n+1}=(u_{n+1}, v_{n+1}, \rho_{n+1}, \psi_{n+1},
\phi_{n+1})$ as solution of the ordinary differential equations
obtained by substituting $u$, $v$, $\rho$, $\psi$, $\phi$ in the
right hand side of the evolution system (\ref{1.14}) to (\ref{18}).
\end{enumerate}
\indent It is very important to notice that, for every $n$ the
initial data for the ordinary differential equations are the same
initial data $u_{0}$, $v_{0}$,  $\rho_{0}$,
 $\psi_{0}$ and   $\phi_{0}$. We obtain through this way a sequence
 $S_{n}=(u_{n}, v_{n}, \rho_{n}, \psi_{n},
\phi_{n})$, $n\in \mathbb{N}$ defined in a maximal interval $[0,
T_{n}[$, $T_{n}$.
\subsection{Boundedness of the iterated sequence}
\begin{proposition}\label{prop1}
There exits $T>0$, T independent on $n$, such that the iterated
sequence $S_{n}=(u_{n}, v_{n}, \rho_{n}, \psi_{n}, \phi_{n})$ is
defined and uniformly bounded over $[0, T_{n}[$, $T_{n}$.
\end{proposition}
\begin{prev}Let $N\in \mathbb{N}$, $N>1$, be an integer. Suppose
that we have, for $n\leq N-1 $, the inequalities
\begin{equation}\label{22}
| u_{n}-u_{0}|\leq C_{1}, \;    |v_{n}-v_{0}| \leq C_{2}, \;
|\rho_{n}-\rho_{0} |\leq C_{3}, \; | \psi_{n}-\psi_{0} |\leq C_{4},
\; | \phi_{n}-\phi_{0} |\leq C_{5}
\end{equation}
 where $C_{i}>0, i=1,...,5$ are given constants. We are going to
 prove that one can choose the constants $C_{i}$ such that
 (\ref{22}) still holds for $n=N$  on $[0, T[$, $T>0$ sufficiently
 small.
\indent Integrating over $[0, t]$, $0\leq t \leq T$, the ordinary
differential equations satisfied by : $u_{N}$, $v_{N}$,  $\rho_{N}$,
 $\psi_{N}$,   $\phi_{N}$ yields :
\begin{equation}\label{23}
|u_{N}-u_{0} |\leq B_{1}t, \; | v_{N}-v_{0} |\leq B_{2}t, \; |
\rho_{N}-\rho_{0} |\leq B_{3}t, \; | \psi_{N}-\psi_{0} |\leq
B_{4}t\leq, \; |\phi_{N}-\phi_{0}|\leq B_{5}t
\end{equation}
where $B_{i}>0, i=1,...5$ are constants depending only on the
constant $C_{i}$. If we choose $T>0$ such that $B_{i}T<C_{i}$,
$i=1,...,5$. Hence for $n=N$, the iterated sequence $(S_{n})$ is
defined and uniformly bounded over $[0, T[$.
\end{prev}
\subsection{Local existence and uniqueness of solution}
\begin{theorem}\label{theo1}
The initial value problem for the Einstein-Scalar Field system on
Bianchi type I space-time has a unique local solution.
\end{theorem}
\begin{prev}We are going to prove that the iterated sequence
$(S_{n})$ converges uniformly on each bounded interval $[0,
\zeta]\subset[0, T[$, $\zeta$, towards a solution $S=(u, v, \rho,
\psi, \phi)$ of the evolution system. For this purpose, we study the
difference $S_{n+1}-S_{n}$. But given the evolution equation
(\ref{1.14}) to  (\ref{18}) in $\phi$ and $\psi$, we will deal with
the difference :
$$\sqrt{2\psi_{n}}-\sqrt{2\psi_{n+1}}=\frac{2(\psi_{n}-\psi_{n+1})}{\sqrt{2\psi_{n}}-\sqrt{2\psi_{n+1}}}.$$
We then need to show first of all that the sequence
$\left(\frac{1}{\sqrt{2\psi_{n}}}\right)$ is uniformly bounded.
\\\indent$\bullet$
By (\ref{18}), the iterated equation providing $\psi_{n+1}$  writes
:
\begin{equation}\label{24}
\dot{\psi}_{n+1}=-2(u_{n}+2v_{n})\rho_{n}
\end{equation}
but by proposition\,\ref{prop1}, there exists a constant $C>0$ such
that we have over $[0, T[$ :
\begin{equation*}
|-2(u_{n}+2v_{n})\psi_{n}|\leq C;
\end{equation*}
(\ref{24}) then gives :
\begin{equation*}
\frac{d\psi_{n+1}}{dt}\geq -C.
\end{equation*}
and integrating over $[0,t]$, $0 \leq t\leq T$ yields:
$$\psi_{n+1}\geq \psi_{0}-Ct.$$
Recall that $\psi_{0} >0 $ ; then taking $t$ sufficiently small such
that $Ct \leq \frac{\psi_{0}}{2}$, we have $\psi_{n+1}\geq
\frac{\psi_{0}}{2}$. Then
$$\dfrac{1}{\sqrt{2\psi_{n+1}}}\leq
\dfrac{1}{\sqrt{\psi_{0}}}$$ which shows that
$\frac{1}{\sqrt{2\psi_{n}}}$ is uniformly bounded over $[0,T[$,
$T>0$ small enough.
\\\indent$\bullet$
Taking the difference between two consecutive iterated equations we
deduce from the evolution equations, using $S_{n}(0)=S_{0}$,
$\forall \, n $, that there exists a constant $C_{2}>0$ such that :
\hspace{-2cm}
\begin{eqnarray}
   &&|u_{n+1}(t)-u_{n}(t)|+|v_{n+1}(t)-v_{n}(t)|+|\rho_{n+1}(t)-\rho_{n}(t)|+|\psi_{n+1}(t)-\psi_{n}(t)|
   \nonumber \\
   &+&|\phi_{n+1}(t)-\phi_{n}(t)| \leq C_{2}\int^{t}_{0}(|u_{n}(s)-u_{n-1}(s)|+|v_{n}(s)-v_{n-1}(s)|+|\rho_{n}(s)-\rho_{n-1}(s)|\nonumber\\
   &+&|\psi_{n}(s)-\psi_{n-1}(s)|+|\phi_{n}(s)-\phi_{n-1}(s)|)ds
   \label{25}
\end{eqnarray}

For the same reasons we have :
\begin{eqnarray}
   &&\left|\frac{du_{n+1}}{dt}(t)-\frac{du_{n}}{dt}(t)\right|+\left|\frac{dv_{n+1}}{dt}(t)-\frac{dv_{n}}{dt}(t)\right|+
   \left|\frac{d\rho_{n+1}}{dt}(t)-\frac{d\rho_{n}}{dt}(t)\right|+\left|\frac{d\psi_{n+1}}{dt}(t)-\frac{d\psi_{n}}{dt}(t)\right|
   \nonumber \\
   &+&\left|\frac{d\phi_{n+1}}{dt}(t)-\frac{d\phi_{n}}{dt}(t)\right|\leq C_{3}(|u_{n}(t)-u_{n-1}(t)|+|v_{n}(t)-v_{n-1}(t)|+|\rho_{n}(t)-\rho_{n-1}(t)|\nonumber\\
   &+&|\psi_{n}(t)-\psi_{n-1}(t)|+|\phi_{n}(t)-\phi_{n-1}(t)|)ds.  \label{26}
\end{eqnarray}
For $n \in \mathbb{N} $, we set :
\begin{equation}\label{27}
\beta_{n}(t)=|u_{n+1}(t)-u_{n}(t)|+|v_{n+1}(t)-v_{n}(t)|+|\rho_{n+1}(t)-\rho_{n}(t)|+|\psi_{n+1}(t)-\psi_{n}(t)|+
|\phi_{n+1}(t)-\phi_{n}(t)|,
\end{equation}
(\ref{25}) and (\ref{27}) give :
\begin{equation}\label{28}
\beta_{n}(t)\leq C_{2}\int^{t}_{0}\beta_{n-1}(t).
\end{equation}
By induction on $n\geq2$, we obtain, from (\ref{28}) :
\begin{equation}\label{29}
|\beta_{n}(t)|\leq \|\beta_{2}\|\dfrac{(C_{2}t)^{n-2}}{(n-2)!}\leq
\|\beta_{2}\|\dfrac{(C_{2}\zeta)^{n-2}}{(n-2)!}
\end{equation}
for $0\leq t \leq  \zeta $ and $0< \zeta < T .$ But the series
$\sum\limits_{n=0}^{+\infty}\dfrac{C^{n}}{n!}$ converges. Hence we
obtain from (\ref{29}) that : $$ \lim\limits_{t \to
+\infty}\sup\limits_{0\leq t \leq \zeta} \beta_{n}(t)=0 .$$
According to definition (\ref{28}) of $\beta_{n}$, we conclude that
every sequence $u_{n}$, $v_{n}$, $\rho_{n}$, $\psi_{n}$ and
$\phi_{n}$ converges uniformly on every interval $[0, \zeta]$,
$0<\zeta<T$ and we denote the different limits by $u$, $v$, $\rho$,
$\psi$ and $\phi$ are continuous functions of $t$.\\
Now from the inequality (\ref{26}), we conclude similarly that the
sequences of derivatives $\left(\frac{du_{n}}{dt}\right)$,
$\left(\frac{dv_{n}}{dt}\right)$,
$\left(\frac{d\rho_{n}}{dt}\right)$,
$\left(\frac{d\psi_{n}}{dt}\right)$,
$\left(\frac{d\phi_{n}}{dt}\right)$ converge uniformly on $[0,
\zeta]$, $0<\zeta < T$. In this conditions, the functions $u$, $v$,
$\rho$, $\psi$ and $\phi$ are of class $\mathcal{C}^{1}$ on $[0,T[$.
Hence $S=(u, v, \rho, \psi, \phi)$ is a local solution of the system
(\ref{1.14}) to (\ref{18}).\\
We now prove that the solution is unique. Consider two solutions
$S_{1}$ and $S_{2}$ of the same initial values problem. Define
$\beta(t)=|S_{1}-S_{2}|$ with $\beta(0)=0$. Since the functions $u$,
$v$, $\rho$, $\psi$ and $\phi$ are bounded on $[0, \zeta]$,
$0<\zeta<T$, there exists a constant $C>0$ such that :
\begin{equation*}
\beta(t)\leq C\int^{t}_{0}\beta(s)ds.
\end{equation*}
By Gronwall lemma, we obtain $\beta(t)=0$ since $\beta(0)=0$,
$S_{1}=S_{2}$ and the local solution is unique. This completes the
proof of proposition\,\ref{prop1}.
\end{prev}
\section{Global existence of solutions}
What we want to know now is whether,  the solution found previously
is global. So, by always following the standard theory on the first
order differential systems, to show that the solution is global, it
will be enough if we prove that $u$, $v$, $\rho$, $\psi$ and $\phi$
remain uniformly bounded.
\begin{remark}
If equations (\ref{1.7}), (\ref{1.8}) admits a global solution
$(a,b)$ defined on $[0, +\infty[$, then $a$ and $b$ will be of class
$\mathcal{C}^{2}$ on $[0, +\infty[$ and hence $u$, $v$, $\rho$,
$\psi$ and $\phi$ are of class $\mathcal{C}^{1}$ on $[0, +\infty[$.
Inversely, if system (\ref{1.14}), (\ref{1.15}), (\ref{1.16}),
(\ref{17}), (\ref{18}) admits a global solution $(u, v, \rho, \psi,
\phi)$ on $[0, +\infty[$, then in particular $u$ and $v$ will be of
class $\mathcal{C}^{1}$ on $[0, +\infty[$ and accordingly the system
(\ref{1.7}), (\ref{1.8}) will admit a global solution of class
$\mathcal{C}^{2}$ on $[0, +\infty[$
\end{remark}
\begin{theorem}\label{theo2}If $\Lambda \geq 0$ and $\dot{b}_{0}>0$, then the
Einstein-Scalar Field system on Bianchi type I space-time has a
global solution.
\end{theorem}
\begin{prev}
Following the standard theory of the first order differential
systems, it will be enough if we prove that every solution of the
Cauchy problem is uniformly bounded. Suppose that $\Lambda\geq0$ and
$\dot{b}_{0}$. \\
$\bullet$ Firstly, the constraint (\ref{19}) implies
$v(v+2u)=\Lambda+8\pi\rho + 8\pi\psi > 0 $ since $8\pi\psi > 0$.
This show that $v$ never vanishes and has the same sign as $v+2u$.
Since, $v$ is continuous and $v(0)= \frac{\dot{b}_{0}}{b_{0}}>0$,
these imply that $v > 0 $ ; the we also have :
$$   v+2u >0 .$$
According to (\ref{17})
$$ \dot{\rho}=-\frac{4}{3}(u+2v)\rho=-\frac{4}{3}\left(\frac{1}{2}(v+2u)+\frac{3}{2}v \right)\rho < 0    $$
since $v+2u >0$, $v >0$ and $\rho >0$, therefore $\rho$ is a
decreasing function on $[0, +\infty[$; it follows that :
$$ 0 < \rho \leq \rho_{0}.$$
$\bullet$ Next according to (\ref{18}),
\begin{eqnarray*}
  \dot{\psi} &=&-2(u+2v)\psi  \\
             &=& -2\left(\frac{1}{2}(v+2u)\psi+\frac{3}{2}u\psi\right) \\
             &=&-(v+2u)\psi-3v\psi < 0,
\end{eqnarray*}
since $v+2u >0$, $v >0 $ and $\psi >0 $ ; it follows that $\psi$ is
a decreasing function on $[0, +\infty[$, hence $$ 0 < \psi \leq
\psi_{0} .$$
 $\bullet$ Finally, let us show that $u$ and $v$ are uniformly
 bounded. Setting $H=u+2v$, it will be enough  to show that H is
 uniformly bounded on $[0, +\infty[$.
\\\indent
We notice that :
$$ H=u+2v=\frac{1}{2}(v+2u)+\frac{3}{2}v \geq \frac{3}{2}v>0  ,$$
so (\ref{1.14}) and  (\ref{1.15}) yield :
\begin{eqnarray}
 \dot{H}  &=& \dot{u}+2\dot{v}  \nonumber  \\
   &=&  2\Lambda-u^{2}-3v^{2}-2uv-8\pi\psi    .\label{30}
\end{eqnarray}
We then have :
\begin{eqnarray*}
 -u^{2}-3v^{2}-2uv  &=& -(u^{2}+v^{2}+4uv)+v^2+2uv \nonumber\\
   &=&-H^{2} + v(v+2u) .\label{31}
\end{eqnarray*}
So (\ref{30}) yields :
\begin{equation}\label{31'}
    \dot{H}= 2\Lambda - H^2 + v(v+2u)-8\pi\psi.
\end{equation}
Let us show that $v(v+2u)$ is bounded above since we have $ v(v+2u)
\geq 0$. Hamiltonian constraint can be written as :
\begin{equation}\label{32}
v(v+2u)-8\pi\psi = \Lambda + 8\pi\rho.
\end{equation}
Setting $ \Lambda_{0}=v(v+2u)-8\pi\psi$ hence $\Lambda_{0}= \Lambda
+ 8\pi\rho $, but we have :
\begin{eqnarray*}
0 < \rho \leq \rho_{0}   &\Longleftrightarrow&  0< 0\pi \rho \leq 8\pi \rho_{0} \\
                          &\Longleftrightarrow& \Lambda\leq\Lambda_{0}\leq \Lambda
+ 8\pi\rho_{0}
\end{eqnarray*}
which means that $\Lambda_{0}$ is bounded. Otherwise we have $0 <
\psi \leq \psi_{0} $ hence $ \Lambda< v(v+2u)\leq \Lambda +
8\pi\rho_{0}+ 8\pi\psi_{0}$ which shows that $v(v+2u)$ is bounded.
(\ref{31'}) implies, since $v(v+2u)-8\pi\psi = \Lambda + 8\pi\rho$,
that
$$ \dot{H} \leq  3\Lambda+ 8\pi\rho_{0} -H^2 .$$
But it is well known that with $C^{2}_{0}= 3\Lambda+ 8\pi\rho_{0} $,
that : $$H(t) \leq W(t),$$ where W satisfies :
\begin{equation}\label{33}
    \left\{
      \begin{array}{ll}
      \dot{W}= C^{2}_{0}-W^2   & \hbox{;} \\
        W(0)=H(0) & \hbox{.}
      \end{array}
    \right.
\end{equation}
We give a general result that is useful here and in what follows.
Consider the Cauchy problem, in which $t_{0} \in \mathbb{R} $ is
given :
\begin{equation}\label{35}
    \left\{
      \begin{array}{ll}
      \dot{y}=K-\alpha^2 y^2   & \hbox{;} \;\;\; (a) \\
      y(t_{0}) \;   \hbox{given}&\;\;\; (b)
      \end{array}
    \right.
\end{equation}
where $K>0$ and $\alpha>0$ are constants ; (\ref{35})(a) is a first
order differential equation of Riccati type, which admits
$y_{0}=\frac{K}{\alpha}$ as an evident solution. It is also well
known that, setting $y=Z+ \frac{K}{\alpha}$ leads to a Bernouilli
equation in Z, which turns out to be a first order linear
differential equation in $\frac{1}{Z}$. Hence, by direct calculation
we obtain :
\begin{equation}\label{36}
    y(t)= \frac{K}{\alpha}\left[1+\frac{2h_{1}(t_{0})}{h_{2}(t_{0})\exp(2\alpha Kt)- h_{1}(t_{0})}\right]
\end{equation}
where : $h_{1}(t_{0})=\alpha y(t_{0})-K$ ; $h_{2}(t_{0})=\alpha
y(t_{0})+K$. (\ref{36}) shows that $y(t)\longrightarrow
\frac{K}{\alpha}$ as $t \longrightarrow + \infty$ thus $y$ is
bounded. Applying this result to (\ref{33}) by setting $K=C_{0}$,
$\alpha =1 $ and $t_{0}=0$, it appears that $W(t)\longrightarrow
C_{0}$ as $t \longrightarrow + \infty$, then W is bounded. The
corresponding reduced expression (\ref{36}) for W ensures that
$W\geq 0 $. Now, since $H(t)\leq W(t)$, $t\geq0$, and since W is
bounded, H is bounded from above. This completes the proof of
Theorem\,\ref{theo2}
\end{prev}

\section*{Acknowledgements}
I acknowledge with thanks the support of the Higher Teacher's
Training College, University of Maroua where this paper was
initiated, prepared and finalized.

\end{document}